\documentclass[12pt]{article}
\usepackage{amsmath,amsfonts,amsthm}

\newcommand{\IN}{{\mathbb{N}}}
\newcommand{\IR}{{\mathbb{R}}}
\newcommand{\F}{{\mathcal{F}}}
\newcommand{\B}{{\mathcal{B}}}

\newtheorem{theorem}{Theorem}[section]
\newtheorem{lemma}{Lemma}[section]
\newtheorem{example}{Example}[section]
\newtheorem{corollary}{Corollary}[section]

\title{Approximating distribution functions by
iterated function systems}

\author{Stefano M. Iacus\footnote{Address for correspondance: Stefano M. Iacus, Dipartimento di Economia e Politica Aziendale, Via Conservatorio 7, I-20122 Milano, Italy. E-mail: {\sl stefano.iacus@unimi.it}} \quad Davide La Torre\\
Universit\`a degli Studi di Milano\\
Dipartimento di Economia Politica e Aziendale\\
Via Conservatorio 7 - I-20122 Milano}
\date{}

\begin{document}
\maketitle

\begin{abstract}
In this paper an iterated function system on the space
of distribution functions is built. The inverse problem is introduced and studied by convex optimization problems. Some applications
of this method to approximation of distribution functions and to estimation theory are given.

\end{abstract}

\noindent
{\bf keywords} : iterated function systems, optimization, distribution function estimation.


\section{Introduction}

The Iterated Function Systems (IFSs) were born in mid eighties
(Hutchinson 1981, Barnsley and Demko, 1985) as applications of the theory of discrete dynamical
systems and as useful tools to build fractals and other similar
sets.
Some possible applications of IFSs can be found in image
processing theory (Forte and Vrscay, 1994), in the theory of stochastic growth
models (Montrucchio and Privileggi, 1999)  and in the theory of random
dynamical systems (Arnold and Crauel 1992, Elton and Piccioni 1992, Kwiecinska and Slomczynski, 2000).
The fundamental result  (Barnsley and Demko, 1985) on which the
IFS method is based is  Banach theorem.

In practical applications a crucial problem is the so-called  {\it inverse
problem}. This can be formulated as follows: given $f$ in some metric space $(S,d)$, find a contraction  $T:S\to S$ that admits a unique fixed point $\tilde f\in
S$ such that $d(f,\tilde f)$ is small enough. In fact if one is
able to solve the inverse problem with arbitrary precision, it is
possible to identify $f$ with the operator $T$ which has
it as fixed point.

The paper is organized as follows: Section \ref{sec2} is devoted to introduce a contractive operator $T$ on the space of distribution functions while, in Section \ref{sec3}, the inverse problem for $T$ 
is studied in details. Section \ref{sec4} is divided into two parts: in the first some examples of inverse problems are analyzed and explicit solutions are given. In the second one, we introduce an estimator of the unknown distribution function based on IFSs.

\section{A contraction on the space of distribution functions}\label{sec2}

Let us denote by the space of distribution functions $F$ on $[0,1]$ by  $\F([0,1])$ and by
$\B([0,1])$  the space of bounded functions on $[0,1]$. Let us further define, for $F,G \in \B([0,1])$,
$d_{sup}(F,G)= \sup_{x\in [0,1]}|F(x)-G(x)|$. So that  $(\F([0,1]),d_{sup})$ is a metric space.
\begin{lemma}
The space $(\F([0,1],d_{sup})$ is a complete metric space.
\end{lemma}
\begin{proof}
Let $F_n$ be a Cauchy's sequence in $\F([0,1])$. Then
 $F_n$ converges to $F$ in $(\B([0,1]),d_{sup})$. Furthermore
it is true that
$F(0)=\lim\limits_{n\to +\infty} F_n(0)=0$ and $F(1)=\lim\limits_{n\to +\infty} F_n(1)=1$
and that if $x_1\ge x_2$ then: 
$$
F(x_1)=\lim\limits_{n\to +\infty} F_n(x_1)\ge \lim\limits_{n\to +\infty} F_n(x_2)=F(x_2).
$$
To prove the right continuity of $F$ we use the uniform convergence of
$F_n$ to $F$, obtaining:
$$
\lim\limits_{x\to a^+} F(x)=\lim\limits_{x\to a^+} \lim\limits_{n\to +\infty} F_n(x)=
\lim\limits_{n\to+\infty}\lim\limits_{x\to a^+} F_n(x)=F(a).
$$
\end{proof}
On $(\F([0,1],d_{sup})$
we define an operator in the following way:
\begin{equation}
TF(x)=p_i F(w_i^{-1}(x)) + \sum_{j=1}^{i-1} p_j + \sum_{j=1}^{i-1}\delta_j\,,
\quad  x\in w_i([a_i,b_i))\,,
\label{func}
\end{equation}
where $F\in \F([0,1])$, $k\in\IN$ is fixed and:
\begin{description}
\item{i)} $w_i:[a_i,b_i)\to [c_i,d_i)=w_i([a_i,b_i))$, $i=1,\ldots, k-1$, 
 $w_k:[a_k,b_k]\to [c_k,d_k]$, with $a_1=c_1=0$ and $b_k=d_k=1$;
\item{iii)} $w_i$, $i=1\ldots k$, are increasing and continuous;
\item{ii)} $\bigcup\limits_{i=1}^k w_i([a_i,b_i))=[0,1)$;
\item{iii)} $p_i\ge 0$, $i=1,\ldots, k$, $\delta_i\ge 0$, $i=1\ldots
k-1$, $\sum\limits_{i=1}^k p_i+\sum\limits_{i=1}^{k-1} \delta_i=1$;
\item{iv)} if $i\not=j$ then $w_i([a_i,b_i))\cap w_j([a_j,b_j))=\emptyset$. 
\end{description}
A similar approach has been discussed in
La Torre and Rocca (1999) but here a more general operator is defined.
In the following we will think
that the maps $w_i$ and the parameters $\delta_j$
are fixed while the parameters $p_i$
have to be chosen. To put in
evidence the dependence of the operator
$T$ on the vector $p=(p_1,\ldots,p_k)$ we will
write $T_p$ instead of $T$.
In many pratical cases $w_i$ are affine maps.
In Corollary \ref{coro} the hypothesis iii) will be weakened to allow more general
functionals.
\begin{theorem}
$T_p$ is an operator from $\F([0,1])$ to itself.
\end{theorem}
\begin{proof}
It is trivial that $T_p F(0)=0$ and $T_p F(1)=1$.
Furthermore if $x_1>x_2$, without
loss of generality, we will consider the two cases:
\begin{description}
\item{i)} $x_1,x_2\in w_i([a_i,b_i))$;
\item{ii)} $x_1\in w_{i+1}([a_{i+1},b_{i+1}))$ and $x_2\in w_i([a_i,b_i))$.
\end{description}
In case $i)$, recalling that $w_i$ are increasing maps, we have:
$$
\begin{aligned}
T_p F(x_1)&=p_i F(w_i^{-1}(x_1)) + \sum_{j=1}^{i-1} p_j +
\sum_{j=1}^{i-1} \delta_j\\
&\ge
p_i F(w_i^{-1}(x_2)) + \sum_{j=1}^{i-1} p_j + \sum_{j=1}^{i-1} \delta_j\\
&=T_p F(x_2)
\end{aligned}
$$
In case $ii)$ we obtain:
$$
T_p F(x_1)-T_p F(x_2)=p_i + \delta_{i-1}+
p_{i+1}F(w_{i+1}^{-1}(x_1))- p_i F(w_i^{-1}(x_2))=
$$
$$
=p_i (1 - F(w_i^{-1}(x_2))) + p_{i+1}F(w_{i+1}^{-1}(x_1))+\delta_{i-1}\ge 0
$$
since $p_i\ge 0$, $\delta_i\ge 0$ and $0\le f(y)\le 1$.
Finally, one can prove without difficulties the right continuity of $T_pf$.
\end{proof}
The following corollary to the previous result will be useful for the applications in Section \ref{sec4}.
\begin{corollary}\label{coro}
Suppose that $w_i:[a_i,b_i) \to [a_i,b_i)$, $w_i(x)=x$, $p_i=p$, $\delta_i\geq -p$, $i=1,\ldots,k$. Then $T_p: \F([0,1]) \to \F([0,1])$.
\end{corollary}
\begin{proof}
Looking at the proof of previous theorem one sees that it is only necessary to prove that
$T_p$ is a non decreasing function. Case i) is analogous whistl for case ii), chosing $x_1>x_2$ we have:
$$
\begin{aligned}
T_p F(x_1) - T_p F(x_2) &= p (1-F(x_2))+p F(x_1)+\delta_i\\
&=p (F(x_1)-F(x_2)) + p + \delta_i \geq 0 
\end{aligned}
$$
\end{proof}
\begin{theorem}
If $c=\max\limits_{i=1,\ldots, k}p_i<1$, then $T_p$ is a
contraction on $(\F([0,1]),d_{sup})$ with
contractivity constant $c$.
\end{theorem}
\begin{proof}
Let $F,G\in (\F([0,1]),d_{sup})$ and let it be $x\in w_i([a_i,b_i))$. We have
$$
|T_p F(x)-T_p G(x)|\le 
p_i\left| F(w_i^{-1}(x)) - G(w_i^{-1}(x))\right|\le c\,
d_{sup}(F,G)\,.
$$
This implies $d_{\infty}(T_p F,T_p G)\le c \,d_{\infty}(F,G)$. 
\end{proof}
The following theorem states that small perturbations of the parameters $p_i$ produce small variations on the fixed point of the operator.
\begin{theorem}
Let $p$, $p^* \in \IR^k$ such that $T_p F_1 =F_1$ and $T_{p^*} F_2 = F_2$. Then
$$d_\infty(F_1,F_2)\leq \frac{1}{1-c}\sum\limits_{j=1}^k  \left| p_j - p_j^*\right| $$
where $c$ is the contractivity constant of $T_p$.
\end{theorem}
\begin{proof}
In fact we have
$$
\begin{aligned}
d_\infty(F_1,F_2) &= d_\infty(T_p F_1,T_p F_2) \\
&=
\max_{i=1,\ldots,k} \left\{
\left|p_i F_1(w_i^{-1}(x)) + \sum\limits_{j=1}^{i-1} p_j -
p_i^* F_2(w_i^{-1}(x)) + \sum\limits_{j=1}^{i-1} p_j^*\right|
\right\}\\
&\leq\sum\limits_{i=1}^k |p_i -p_i^*| + c\, d_\infty(F_1,F_2)\,,
\end{aligned}
$$
since
$$
\begin{aligned}
&\left|p_i F_1(w_i^{-1}(x)) + \sum\limits_{j=1}^{i-1} p_j -
p_i^* F_2(w_i^{-1}(x)) + \sum\limits_{j=1}^{i-1} p_j^*\right|\\
&\leq
\sum\limits_{j=1}^{i-1} |p_j -p_j^*| +
|p_i F_1(w_i^{-1}(x)) - p_i F_2(w_i^{-1}(x))|\\
&\phantom{ss}+|p_i F_2(w_i^{-1}(x)) - p_i^* F_2(w_i^{-1}(x))|\\
&\leq
\sum\limits_{j=1}^{i-1} |p_j -p_j^*|  + p_i d_\infty(F_1,F_2) + 
|p_i-p_i^*|\\
&\leq c\,d_\infty(F_1,F_2) + \sum_{j=1}^k |p_j-p_j^*|\,.
\end{aligned}
$$
\end{proof}

\section{The inverse problem as a convex constrained optimization problem}\label{sec3}
Choose $F\in (\F([0,1]),d_{sup})$. The aim of solving the inverse problem 
is to find a contractive map $T:\F([0,1])\to \F([0,1])$ which has a
fixed point ``near"
to $F$.
In fact if it is possible to solve the inverse problem with an arbitrary
precision one can identify the operator $T$ with its fixed point.
With a fixed system of maps $w_i$ and parameters $\delta_j$, the inverse
problem can be solved, if it is possible, by using the parameters $p_i$.
These have to be choose in the following convex set:
$$
C=\left\{ p\in\IR^k: p_i\ge 0, i=1,\ldots, k,
\sum_{i=1}^k p_i=1-\sum_{i=1}^{k-1} \delta_i \right\},
$$
We have the following result.
\begin{theorem}
\label{approx}
Choose $\epsilon>0$ and $p\in C$ such that  $p_i\cdot p_j>0$ for some $i\neq j$.  
If $d_{sup}(T_p F, F)\le\epsilon$, then: 
$$
d_{sup}(F,\tilde F)\le \frac{\epsilon}{1-c},
$$
where $\tilde F$ is the fixed point of $T_p$ on
$\F([0,1])$ and $c=\max\limits_{i=1,\ldots, k} p_i$  is
the contractivity constant of $T_p$.
\end{theorem}

\begin{proof}
The assumptions imply $c<1$. So we have:
$$
d_{sup}(F,\tilde F)\le d_{sup}(F,T_pF) +
d_{sup}(T_pF,T_p\tilde F)\le \epsilon +
c\, d_{sup}(F,\tilde F)
$$
and  so we get the thesis. 
\end{proof}

If we wish to find an approximated solution of the inverse problem,
we have to solve the following constrained optimization problem:
$$
\min_{p\in C} d_{sup}(T_p F,F) \leqno{({\bf P})}
$$

It is clear that the ideal solution of ({\bf P}) consists of finding a $p^*\in C$ such
that $d_{sup}(T_{p^*} F,F)=0$. In fact this means that, given a distribution
function $F$, we have found a contractive map $T_p$ which has
exactly $F$ as fixed point. Indeed the use of Theorem \ref{approx}
gives us only an approximation of $F$. This can be improved,
once fixed the maps $w_i$,
increasing the number of parameters $p_i$.

The following result 
proves the convexity of the function 
$D(p)=d_{sup}(T_pF,F)$, $p\in\IR^k$.
\begin{theorem}
The function $D(p):\IR^k\to\IR$ is convex.
\end{theorem}
\begin{proof}
If we choose $p_1, p_2\in\IR^k$ and $\lambda\in [0,1]$ then: 
$$
D(\lambda p_1 + (1-\lambda) p_2) = \sup_{x\in [0,1]} |T_{\lambda p_1 + (1-\lambda) p_2} F(x)-F(x)|
\le
$$
$$
\lambda \sup_{x\in [0,1]}|T_{p_1} F(x)-F(x)| + 
(1-\lambda) \sup_{x\in [0,1]}|T_{p_2} F(x)-F(x)|= \lambda D(p_1) + (1-\lambda) D(p_2).
$$
\end{proof}   
Hence for solving  problem ({\bf P}) one can recall classical results about convex 
programming problems (see for instance Rockafellar and Wets, 1998). A necessary and
sufficient condition for $p^*\in C$
to be a solution of ({\bf P})
can be given by Kuhn-Tucker conditions.

\section{Inverse problem for distribution functions and applications}\label{sec4}
In this section we consider different problems. We show that for a particular class of distribution functions the inverse problem can be solved exactly without solving any optimization problem.
Then we discuss two ways of construct IFS to approximate a distribution function $F$
with a finite number of parameters $p_i$ and maps $w_i$.

As is usual in statistical applications, given a sample of $n$ independent and identically distributed  observations, $(x_1,x_2, \ldots,x_n)$, drawn from an unknown distribution function $F$, one can easily contruct the empirical distribution function (e.d.f.) $\hat F_n$ that reads
$$
\hat F_n(x) = \frac1{n} \sum\limits_{i=1}^n \chi_{(-\infty,x]}(x_i),\quad x\in\IR\,,
$$
where $\chi_A$ is the indicator function of the set $A$.
Asymptotical properties of optimality of $\hat F_n$ as estimator of the unknown $F$ when $n$ goes to infinity are well known and studied (Millar 1979 and 1983).
This function has an IFS representation that is exact and can be found without solving any optimization problem. We assume that the 
 the $x_i$ in the sample are all different (this assumption is natural if $F$ is a continuous distribution function).
Let $w_i(x) : [x_{i-1},x_i)\to [x_{i-1},x_i)$, 
when $i=1\ldots n$ and $w_1(x):[0,x_1)\to [0,x_1)$,
$w_{n+1}(x):[x_n,x_{n+1}]\to [x_n,x_{n+1}]$,with $x_0=$ and $x_{n+1}=1$.
Assume also that every map is of the form $w_i(x)=x$. If
we choose $p_i=\frac{1}{n}$, $i=2\ldots n+1$,
$p_1=0$ and 
$$\delta_1 = \frac{n-1}{n^2},\quad \delta_i=-\frac{1}{n^2}$$
then the following representation holds:
$$T_p \hat F_n(x)=
\begin{cases}
0,& i=1\\
\frac{1}{n} \hat F_n(x) + \frac{n-1}{n^2},&i=2\\
 \frac{1}{n}\, \hat F_n (x) + \frac{i-1}{n} + \frac{n-i+1}{n^2},& i=3,\ldots,n+1.
\end{cases}
$$
when $x\in [x_{i-1}, x_i)$.
Furthermore, from Corollary \ref{coro} we are guaranteed that 
$$\lim_{s\to\infty} d_\infty(T_p^{(s)} u, \hat F_n) \to 0,\quad \forall\, u\in\mathcal F[0,1].
$$
Note that from the point of view of applications, constructing the e.d.f. or iterate the IFS with the given maps is exactly equivalent if one start, for example, with a uniform distribution on [0,1] in the first iteration.
So this is just a case when we can present an IFS that gives exact result for this particular class of distribution functions.

What follows, on the contrary, is more attractive from the point of view of applications. Suppose that one knows the distribution function $F$ and wants to construct the IFS which has $F$ as fixed point. In general one has to provide an infinite set of affine maps $\{w_i, i\in \IN\}$ and solve an extremal problem to find the corresponding sequence of weigths $p_i$, $i\in\IN$. This problem has not a general solution but at the same time the solution in terms of a finite, possibily few, number of maps and weigths is crucial in applications like image compression and trasmission.

The idea is the following: one can think at $n$ points $(x_1,x_2,\ldots,x_n)$ as they were drawn   from the distribution function $F$, and use the same maps $w_i$ of the e.d.f. $\hat F_n$, then instead of using the $p_i$ equal to $1/n$ one solve the extremal problem as it is usual in IFS application.
The corresponding IFS should have a fixed point that is a ``good" approximation on $F$. So it is sufficient to store the simulated data and the weights instead of $F$ itself. 

We take the functional $T_p F$ with the particular choice of $\delta_i = 0$. This choice is in principle not necessary but simplifies the solution of the problem. We simulated $n$ i.i.d. observations from the distribution function $F$ and we use the maps of the e.d.f. above.

We now try to solve the extremal problem 
$$
\min d_\infty(T_p F,F) 
$$
under the constrain $\sum\limits_{i=1}^n p_i=1$, $p_i\geq 0$, $i=1,\ldots,n$ (with some $p_i >0$).
The optimal solution will be $\{\hat p_i, i=1,\ldots,n\}$ such that $d_\infty(T_{\hat p} F,F)=0$ that it is true in at least one case: if $F$ equals $\hat F_n$ and $p_i = 1/n$. Otherwise we will obtain some positive number. That means that, in principle, in the worst case we can approximate $F$ with its empirical distribution function $\hat F_n$.
But we can generally do better. So let us solve the problem: let us fix $x_0=0$ and $x_{n+1}=1$, then
$$
\begin{aligned}
d_\infty(T_p F, F) &= \sup_{x\in[0,1]} \left |  T_p F(x) - F(x) \right |\\
&=\max_{i=1,\ldots,n+1} \left\{
\sup_{[x_{i-1},x_i)} \left | T_p F(x) - F(x)\right| 
\right\} \\
&=\max_{i=1,\ldots,n+1} \left\{
\sup_{[x_{i-1},x_i)} \left |\sum\limits_{j=1}^{i-1} p_j - (1-p_i)  F(x) \right| 
\right\} \\
&=\max_{i=1,\ldots,n+1} \left\{
 \left | \sum\limits_{j=1}^{i-1} p_j - (1-p_i)  F(x_{i-1}) \right| ,
 \left | \sum\limits_{j=1}^{i-1} p_j - (1-p_i)  F(x_i^-) \right| 
\right\}
\end{aligned}
$$
and the last line is due to the non-decreasing property of $F$.

\begin{example}
Suppose that $F$ is the distribution function of a uniform distribution on [0,1] and suppose that we can only draw one observation from $F$ (or choose a point) $x_1$. The empirical distribution function $\hat F_1(x) = \chi_{\{x_1,1\}}(x)$ is usless if we have in mind to approximate $F$.
Let us use the second technique: fix $w_1 : [0,x_1) \to [0,x_1)$ and $w_2 : [x_1,1] \to [x_1,1]$. 
We try to solve the above extremal problem.
$$
\begin{aligned}
d_\infty(T_p F, F) &=\max_{i=1,\ldots,n+1} \left\{
 \left | \sum\limits_{j=1}^{i-1} p_j - (1-p_i)  F(x_i) \right| ,
 \left | \sum\limits_{j=1}^{i-1} p_j - (1-p_i)  F(x_{i+1}) \right| 
\right\}\\
&= \max
\biggl \{ |-(1-p_1)\cdot 0 |, |-(1-p_1)x_1|,\\
&\phantom{=}|p_1-(1-p_2)x_1|,|p_1-(1-p_2)|
\biggr\}\\
&=\max\{0, x_1(1-p_1), p_1(1-x_1),0 \},\qquad x_1\in(0,1), \, p_1+p_2=1\,.
\end{aligned}
$$
now minimazing with respect to $p_1$ and $p_2$ under the constrain $p_1+p_2=1$ one obtains
simply $p_1=x_1$.
The resulting functional will be
$$
T_{x_1} u(x) =
\begin{cases}
x_1 \, u(x),& x\in [0,x_1)\\
(1-x_1)\, u(x) + x_1,& x\in [x_1,1]
\end{cases}
$$
and it is clear that $T_{x_1} u(x) = T_{x_1} x$, for one iteration only, is closer than $\hat F_1$ to $F(x) = x$ and the approximation is better and better as $x_1\to 0$ or $x_1\to 1$.
\end{example}

We propose now a more efficient method to approximate $F$ when $F$ is not to be estimated.
We have already mentioned that the e.d.f. is the better estimator of an unknown distribution function $F$, so one can think to sample $n$ points from $F$ and use their values to approximate $F$ by $\hat F_n$. As $n\to\infty$, the statistical literature assures almost sure convergence of $\hat F_n(x)$ to $F(x)$ for every $x$.
We also have shown the exact IFS representation of  $\hat F_n$.
But this method is not efficient.
On the contrary, suppose that $F$ is a continuous distribution function. As we know $F$, we can think to approximate it by means of continuous functions instead of simple functions like $\hat F_n$.
Choose $n$ points  $(x_1,\ldots,x_n)$ and assume that $x_0=0$ and $x_{n+1}=1$. One can costruct the following functional
$$
T_F u(x) = (F(x_i)-F(x_{i-1})) u\left(\frac{x-x_{i-1}}{x_i-x_{i-1}}\right) + F(x_{i-1}),\quad x\in[x_{i-1},x_i),
$$
$ i=1,\ldots,n+1$.
Notice that $T_F$ is a particular case of \eqref{func} where
$p_i = F(x_i) - F(x_{i-1})$, $\delta_i=0$ and $w_i(x) : [0,1)\to [x_{i-1},x_i)$.
This is a contraction and, at each iteration, $T_F$ passes exactly through the points $F(x_i)$.
It is almost evident that, when $n$ increases the fixed point of the above functional will be ``close" to $F$. So again, instead of sending an infinite set of weigths and maps, one can send $n$
points and the values of $F$ evaluated at these points. All in summary, only $2\cdot n$ informations should be sent to reconstruct $F$.

For $n$ small, the choice of a good grid of point is critical. So one question arises: how to choose the $n$ points ? One can proceed case by case but as $F$ is a distribution function one can use its properties. We propose the following solution: take $n$ points $(u_1, u_2,\ldots, u_n)$ equally spaced $[0,1]$ and define $x_i = F^{-1}(u_i)$, $i = 1,\ldots, n$. The points $x_i$ are just the quantiles of $F$. In this way, it is assured that the profile of $F$ is followed as smooth as possible. In fact, if two quantiles $x_i$ and $x_{i+1}$ are relatively distant each other, than $F$ is slowly increasing in the interval $(x_i,x_{i+1})$ and viceversa.
This method is more efficient than simply taking equally spaced points on $[0,1]$.
If this method of choosing the points is used, the the functional simply appears as
$$
T_F u(x) = \frac{1}{n} u\left(\frac{x-x_{i-1}}{x_i-x_{i-1}}\right) + \frac{i-1}{n},\quad x\in[x_{i-1},x_i),\, i=1,\ldots,n+1\,.
$$
And this suggests an empirical estimator of $F$. If $q_i$, $i=1,2,\ldots,k$, $k<n$, are the  empirical quantiles of the sample $(x_1,x_2,\ldots,x_n)$ of order $\frac{i}{k}$, then an estimator of the unknown distribution function $F$ should be written as
$$
\tilde F_{(k)} u(x) = \frac{1}{k} u\left(\frac{x-q_{i-1}}{q_i-q_{i-1}}\right) + \frac{i}{k},\quad x\in[q_{i-1},q_i),
$$
$ i=1,\ldots,k$, with $q_0=$ and $q_{k+1}=1$
As $n$ and  $k=k(n)$ go to infinity $\hat F_{(k)}$ converges to $F$.  
Relative efficiency of $\tilde F_{(k)}$ with respect to $\hat F_n$ is investigated via simulations. The results are reported in Table \ref{tab1} for differently shaped distribution functions and sample sizes. What emerges is that $\tilde F_{(k)}$ is equivalent to the e.d.f. in the sense of the sup-norm.



\section*{References}
\begin{itemize}
\item[] Arnold, L., Crauel H. (1992),  Iterated function systems and
multiplicative ergodic theory, {\sl Diffusion processes and related
problems in analysis} (Pinsky M.A., Vihstutz V. eds.), vol. II,
Stochastics flows, Progress in probability 27, Birk\"auser, Boston, 283-305.
 
\item[] Barnsley, M.F., Demko, S. (1985),  Iterated function systems and the
global construction of fractals, {\sl Proc. Roy. Soc. London, Ser A}, {\bf 399}, 243-275.

\item[] Cabrelli, C.A., Forte, B., Molter, U.M., Vrscay, E.R. (1992),  Iterated
fuzzy set systems: a new approach to the inverse problem for
fractals and other sets, {\sl J. Math. Anal. Appl.}, {\bf 171},
79-100.
 
\item[] Elton, J.H., Piccioni, M. (1992),  Iterated function systems arising
from recursive estimation problems, {\sl Probab. Theory Related Fields},
{\bf 91}, 103-114.

\item[] Forte, B., Vrscay, E.R. (1994),  Solving the inverse problem for function/image approximation using iterated function systems, I. Theoretical basis, {\sl Fractal}, {\bf 2}, 3, 325-334.

\item[] Hutchinson, J., (1981),  Fractals and self-similarity, {\sl Indiana Univ.
J. Math.}, {\bf 30}, 5, 713-747.

\item[] Kwiecinska, A.A., Slomczynski, W. (2000),   Random dynamical systems
arising from iterated function systems with place-dependent
probabilities, {\sl Stat. and Prob. Letters}, {\bf 50},  401-407.

\item[] La Torre, D., Rocca, M. (1999), Iterated function systems and optimization problems, {\sl Rend. Sem. Mat. Messina SER II}, {\bf 6}, 21, 165-173.

\item[] Millar, P.W. (1979), Asymptotic minimax theorems
        for the sample distribution function, {\it Z. 
        Wahrscheinlichkeitstheorie Verw. Geb.}, {\bf 48}, 233-252.

\item[] Millar, P.W. (1983), The minimax principle in 
        asymptotic statistical theory, {\it Lecture Notes in Math.},
        {\bf 976}, 76-265.

\item[] Montrucchio, L., Privileggi, F. (1999),  Fractal steady states in
stochastic optimal control models, {\sl Ann. Op. Res.},
{\bf 88}, 183-197.

\item[] Rockafellar, R.T., Wets, R.J-B. (1998), {\sl Variational analysis},
Springer.

\end{itemize}
\newpage

\begin{table}
\begin{center}
\begin{tabular}{c|c|c|c|c}
number of points& $d_\infty\left(\tilde F_{(k)}^{(4)} u, F\right)$ &  $d_\infty\left(\hat F_n, F\right)$ &$\frac{(a)}{(b)}\cdot 100\%$ &distribution $F$\\
drawn from $F$ & (a) & (b) & &\\
\hline
10 &  0.20232 & 0.24103 & 83.94\% & Beta(2,2)\\
50 &  0.09376 & 0.10241 & 91.56\% &Beta(2,2)\\
100 & 0.06989  & 0.07131  & 98.01\%& Beta(2,2)\\
500 & 0.02884 & 0.02917  & 98.87\% &Beta(2,2)\\
1000 & 0.02475 & 0.02506 & 98.78\%&Beta(2,2)\\
\hline
10 & 0.18747 & 0.19472 & 96.27\% &  Beta(3,3)\\
50 & 0.09945 & 0.09777 & 101.72\% & Beta(3,3)\\
100 &    0.07103 & 0.07521 & 94.44\%& Beta(3,3)\\
500 &  0.03077 & 0.03061& 100.52\%&Beta(3,3)\\
1000 &   0.01993 & 0.02018 & 98.74\%& Beta(3,3)\\
\hline
10 & 0.20842 & 0.22220 & 93.80\% & Beta(5,3)\\
50& 0.10615 & 0.10517 & 100.93\%& Beta(5,3)\\
100 &   0.06881 & 0.07096& 96.96\%& Beta(5,3)\\
500 &  0.02959  & 0.02971  & 99.60\%& Beta(5,3)\\
1000 & 0.02176 & 0.02194 & 99.17\%&Beta(5,3)\\
\hline
10 & 0.23054 & 0.23301 & 98.94\%  &Beta(3,5)\\
50& 0.08993 & 0.089347 & 100.66\% & Beta(3,5)\\
100 & 0.06541 & 0.06515 & 100.40\% &Beta(3,5)\\
500 & 0.02978 & 0.03015 & 98.80\% & Beta(3,5)\\
1000 & 0.01978 & 0.02003 & 98.77\%&Beta(3,5)\\
\hline
10 & 0.20522 & 0.24492 & 83.79\%&Beta(1,1)\\
50& 0.10456 & 0.11990 & 87.20\% &Beta(1,1)\\ 
100 & 0.07621 & 0.08124 & 93.81\%&Beta(1,1)\\
500 & 0.02938 & 0.02974 & 98.77\%& Beta(1,1)\\
1000 &   0.02382  & 0.02428  &98.09\%  & Beta(1,1)\\
\end{tabular}
\caption{Simultation results. Values are the arithmetic means over 30 trials. The functional is iterated 4 times starting with the uniform distribution on [0,1] as initial point. Functions are evaluated at 20 equally spaced points on [0,1]. The proposed estimator can be said to be almost equivalent as the best estimator $\hat F_n$.}
\end{center}
\label{tab1}
\end{table}


\end{document}